\documentclass[10pt]{amsart}
\pdfoutput=1
\usepackage{amsthm}

\usepackage{amsmath}
\usepackage{wasysym}
\usepackage{pdfpages}

\usepackage[matrix, arrow,all,cmtip,color]{xy}  

\usepackage{hyperref}

\newcommand{\bi}{\begin{enumerate}}
\newcommand{\ei}{\end{enumerate}}

\newcommand{\bp}{\begin{proof}}

\newcommand{\ep}{\end{proof}}
\def\comment#1{}

\def\path{\Omega}

\def\xra{\xrightarrow}

\def\path{\Omega}

\def\eii{\end{itemize}}
\def\bii{\begin{itemize}}
\def\bid{\begin{description}}
\def\eid{\end{description}}

\usepackage{epigraph}
\usepackage{geometry}

\makeatletter
\usepackage{amssymb, stmaryrd}
\makeatother

\def\0{\varnothing}
\def\rtt{\,\rightthreetimes\,}
\def\lra{\longrightarrow}

\def\ilim{{\raisebox{0pt}{$\bigcirc$}} \kern -0.31cm \hbox{$\Yleft$}}

\DeclareMathAlphabet{\mathpzc}{OT1}{pzc}{m}{it}

\def\rtt{\rightthreetimes}
\def\lra{\longrightarrow}

\def\bee{\begin{equation}}
\def\eee{\end{equation}}

\def\ei{\end{enumerate}}
\def\bi{\begin{enumerate}}
\def\eie{\end{enumerate}}
\def\bie{\begin{enumerate}}
\def\ev{\end{verbatim}}
\def\bv{\begin{verbatim}}
\def\bc{\begin{verbatim}}
\def\ec{\end{verbatim}}
\title[Lifting property as negation]
{Point-set topology as diagram chasing computations
}
\author[Misha Gavrilovich\\
\ \ \ \\ \ \ \ \  \\ \ \ \ \ \ \ \ \ \ \\
\textbf{{\footnotesize\tiny To}  Grigori Mints Z"L  {\footnotesize\tiny In memoriam}}]{Misha Gavrilovich
\texttt{m\!\!\!mi\!\!\!is\!\!\!sh\!\!\!h\!\!\!ha\!\!\!ap\!\!\!p@\!\!\!@s\!\!\!sd\!\!\!df\!\!\!f.\!\!\!.o\!\!\!or\!\!\!rg\!\!\!g https://mishap.sdf.org}. 
\\{\tiny.} 
\\ 
{{\tiny To}  \textsf{Grigori Mints Z"L} \tiny In memoriam}}
\begin{document}\large
\begin{abstract}

We observe that  some natural mathematical definitions  are lifting properties 
relative to simplest counterexamples, namely the definitions of surjectivity and injectivity of maps, 
as well as of being connected, separation axioms $T_0$ and $T_1$ in topology, 
having dense image,  induced (pullback) topology, and  every real-valued function  being bounded (on a connected domain). 

 We also offer a couple of brief speculations on cognitive and AI aspects of this observation,
 particularly that in point-set topology some arguments read as diagram chasing computations
 with finite preorders.

%$\emptyset\lra \{\bullet\} \rtt X\lra Y$ and $\{\bullet,\bullet\}\lra\{\bullet\}\rtt X\lra Y$, and
%$X\lra\{\bullet\}\rtt \{\bullet,\bullet\}\lra\{\bullet\}$ define surjective, injective, and connected, resp.
% We list a few more examples from topology of lifting properties relative to simplest counterexamples
%  -- connected, separation axioms $T_0$, $T_1$, dense image,  induced(pullback) topology,
%   every real-valued function  is bounded (on a connected domain) -- and offer

\end{abstract}

\maketitle 
\setlength{\epigraphwidth}{0.45\textwidth}
\epigraph{\sf{2. Hawk/Goose effect. A baby chick does not have any built-in image of
``deadly hawk" in its head but distinguishes frequent, hence, harmless shapes,
sliding overhead from potentially dangerous ones that appear rarely.
Similarly to ``first", ``frequent" and ``rare" are universal concepts that were
not specifically designed by evolution for distinguishing hawks from geese.
This kind of universality is what, we believe, turns the hidden wheels of the
human thinking machinery.}}{Misha Gromov, {\em Math Currents in the Brain.}
}

\section{Introduction. Structure of the Paper}

This note was written for the {\em De Morgan Gazette} 
[DMG] to show that some natural definitions are 
lifting properties relative to the simplest counterexample, and to suggest a way to ``extract'' these lifting properties 
from the text of the usual definitions and proofs. The exposition is in the form of a story and aims to be self-contained 
and accessible to a first year student    
who has taken some first lectures in naive set theory, topology, and who has heard a definition of a category.
A more sophisiticated reader  may find it more illuminating to recover our formulations herself
from reading either the abstract, or the abstract and the opening sentence of the next section. 
The displayed formulae and Figure~1(a) defining the lifting property provide complete formulations of our theorems 
to such a reader.

%A more mathematically inclined reader
%may find it a nice? and easy? exercise not to read the paper but rather to recover the lifting properties by herself.
%The exercise becomes trivial if one looks the displayed formulae and the definition of the lifting property (Figure~1(a)).
%We end the paper by speculations on cognitive and AI aspects. 

Our approach naturally leads to    a more general observation that in basic
point-set topology, a number of arguments are computations based on symbolic
diagram chasing with finite preorders likely accessible to a theorem prover; 
because of lack of space, 
we discuss it in a separate note [G0].
\pagebreak

%No effort has been made to provide a complete? bibliography; the author shall happily include any references
%sent in the next version (if any). 

\section{\label{surinj} Surjection and injection}

We try to find some ``algebraic"  notation to (re)write the {\em text} of the
definitions of surjectivity and injectivity of a function, as found in any standard
textbook. 
%We want it purely formal and straightforward ---  no words, best few letters,
We want something very straightforward and syntactic --- notation for what we
(actually) say, for the text we write,  and not for its meaning, for who knows
what meaning is anyway? 
%and it should read(decipher) to the actual sentences defining these notations.

%(misha: our notation should translate (unwind) to the text of the definitions as it is found in any standard textbook.)

\bii\item[(*)${}_{\text{words}}$] ``A function $f$ from $X$ to $Y$ is {\em surjective} iff for every element $y$ of $Y$ there is an element $x$ of $X$ 
such that $f(x)=y$."
\eii A function from $X$ to $Y$ is an arrow $X\lra Y$. Grothendieck(?) taught us that a point, say ``$x$ of $X$", is 
(better viewed as) as  $\{\bullet\}$-valued point, that is an arrow $\{\bullet\}\lra X$
from a(the?) set with a unique element; 
similarly ``$y$ of $Y$" we denote by
an arrow $\{\bullet\}\lra Y$. Finally, make dashed the arrows  required to  ``exist". We get the diagram Fig.~1(b) 
without the upper left corner; there ``$\{\}$'' denotes the empty set with no elements listed inside of the brackets.
\bii
\item[(**)${}_{\text{words}}$] ``A function $f$ from $X$ to $Y$ is {\em injective} iff  no pair of different points is
sent to the same point of $Y$." 
\eii% of $X$ such that $f$ sends $a$ and $b$ $f(a)=f(b)$."
``A function $f$ from $X$ to $Y$" is an arrow $X\lra Y$.  ``A pair of points" is a
$\{\bullet,\bullet\}$-valued point, that is an arrow $\{\bullet,\bullet\}\lra X$ from a two element set; 
we ignore ``different" for now. 
``the same point" is an arrow $\{\bullet\}\lra Y$. Represent ``sent to" by an arrow $\{\bullet,\bullet\}\lra \{\bullet\}$. 
What about ``different"? if the points are not ``different", then they are "the same" point, 
that is an arrow $\{\bullet\}\lra X$. 

Now all these arrows combine nicely into diagram Figure~1(c). 
How do we read it? We want this diagram to have the meaning of the sentence
(**)${}_{\text{words}}$ above, so we interpret such diagrams as follows:
%(recall that ``commutative'' in category theory
%means that the composition of the arrows along a directed path depends only on the endpoints of the path): 
\bii\item[$(\rtt)$]
``for every commutative square (of solid arrows) as shown  there is a diagonal (dashed) arrow making the 
total diagram commutative" (see Fig.~$1(a)$) 
\eii
(recall that ``commutative'' in category theory
means that the composition of the arrows along a directed path depends only on the endpoints of the path)

Property $(\rtt)$ has a name and is in fact quite well-known [Qui]. It is called {\em the lifting property}, 
or sometimes {\em orthogonality of morphisms},
and is viewed as the property of the two downward arrows; we denote it by $\rtt$.

Now we rewrite (*)${}_{\text{words}}$ and (**)${}_{\text{words}}$ as:
$$\!\!\mathrm{(*)_{\rtt}}\ \ \ \ \ \ \ \ \ \ \ \ \ \ \ \ \ \ \{\} \lra \{\bullet\} \rtt  X\lra Y$$
$$\mathrm{(**)_{\rtt}}\ \ \ \ \ \ \ \ \ \ \ \ \ \ \ \  \{\bullet,\bullet\}\lra \{\bullet\} \rtt  X\lra Y$$ 

%\bii
%\item[(*)] 
%%\item  \label{sur} $X\lra Y$ is surjective iff\ \
%       $ \{\} \lra \{\bullet\} \rtt  X\lra Y$ %(cf.~Fig.~$\ref{fig1}b$)
%       %\item \label{non-empty} $X$ is non-empty or $X=Y=\emptyset$ iff\ \ 
%       %      $  X\lra Y \rtt  \emptyset\lra \{\bullet\}$ (cf.~Fig.~$\ref{fig1}c$)
%
%\item[(**)]      
%      %\item \label{inj}  $X\lra Y$ is injective iff\ \
%              $  \{\bullet,\bullet\}\lra \{\bullet\} \rtt  X\lra Y$ %(cf.~Fig.~$\ref{fig1}d$)
%
%\eii
%

\def\rrt#1#2#3#4#5#6{\xymatrix{ {#1} \ar[r]^{} \ar@{->}[d]_{#2} & {#4} \ar[d]^{#5} \\ {#3}  \ar[r] \ar@{-->}[ur]^{}& {#6} }}
\begin{figure}
\begin{center}
$\ \ \ (a)\ \xymatrix{ A \ar[r]^{i} \ar@{->}[d]_f & X \ar[d]^g \\ B \ar[r]_-{j} \ar@{-->}[ur]^{{\tilde j}}& Y }$ \ \  
%$\rrt ABXY$\ \ \ \
$(b)\  \rrt  {\{\}}  {.} {\{\bullet\}}  X {\therefore} Y $ \ \ \
$(c)\  \rrt {\{\bullet,\bullet\}} {.} {\{\bullet\}}  X {\therefore} Y $\ \ \ 
$(d)\  \rrt X {\therefore} {Y} {\{x,y\}} {.} {\{x=y\}}$\ \ 
\end{center}
\caption{\label{fig1}
Lifting properties. Dots $\therefore$ indicate free variables, i.e.~a property of
what is being defined.
(a) The definition of a lifting property $f\rtt g$: for each $i:A\lra X$ and $j:B\lra Y$
making the square commutative, i.e.~$f\circ j=i\circ g$, there is a diagonal arrow $\tilde j:B\lra X$ making the total diagram 
$A\xra f B\xra {\tilde j} X\xra g Y, A\xra i X, B\xra j Y$ commutative, i.e.~$f\circ \tilde j=i$ and $\tilde j\circ g=j$.
(b) $X\lra Y$ is surjective 
(c) $X\lra Y$ is injective; $X\lra Y$ is an epicmorphism if we forget %never use
that $\{\bullet\}$ denotes a singleton (rather than an arbitrary object
and thus $\{\bullet,\bullet\}\lra\{\bullet\}$ denotes an arbitrary morphism $Z\sqcup Z\xra{(id,id)} Z$) 
(d) $X\lra Y$ is injective, in the category of Sets; $\pi_0(X)\lra\pi_0(Y)$ is injective, in the category
of topological spaces.
}
\end{figure}
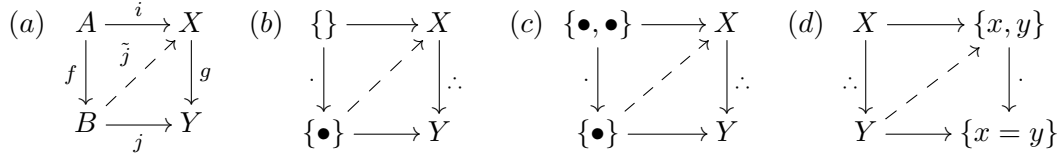

So we rewrote these definitions without any words at all. Our benefits?  
The usual little miracles happen:  

%if we add a top left vertex $\{\}=\emptyset$
%(and arrows $\{\}\lra\{\bullet\}$ and $\{\}\lra X$) 
%to the triangle in Fig.0a (that's Fig.0a'). 

Notation makes apparent a similarity of (*)${}_{\text{words}}$ and (**)${}_{\text{words}}$: they are obtained, in the same purely formal way, from 
the two simplest arrows (maps, morphisms) in the category of Sets. More is true:
it is also apparent that these arrows are the simplest {\em counterexamples} to the properties,
and this suggests that we think of the lifting property as  a category-theoretic (substitute for) negation.
Note also that a non-trivial (=non-isomorphism) morphism never has the lifting property 
relative to itself, which fits with this interpretation. 
%(Confusingly, there are {\em two} negations depending whether the morphism appears on the left or right 
%side of the square.) 

%This sheds some light on a teacher's explanations of surjective and injective:
%every point lifts to its preimage, points do not get glued together. (fixme: do they?)
%\footnote{fixme: I suppose I also want to mention what I described here: 
%The teacher draws two ovals representing sets $X$ and $Y$,
%one above the other, and then an arrow between them going downwards.
%To explain injection, she says "here are two points in $X$ going into one in $Y$,
%and this does not happen"
%(draws two  points inside the upper oval, and a point inside the lower oval $Y$ and then
% joins the upper two to the third one);
% to explain surjection, "here is one point in $Y$ and it must come from somewhere in $X$"
% (draws one point, then a point above it, and an arrow going downwards).
%  
%
%  Fig.1(b-c), Fig.2(top)  seem a fair diagrammatic representation of these explanations, and
%  remarkably like what's actually drawn on the board if
%  you think of a point inside of an oval as an arrow from that point to the oval.
%  (Incidentally, a point refers to an arrow, hence an $X$-valued point of $Y$ is
%  by definition an arrow $X\lra Y$.)
%}

%both (*) and (**) are defined as classes of arrows left orthogonal (left-lifting) to 
%the simplest counterexample to the property. Note that gives some formal sense 
%to a teacher's remark explaining injectivity: the points do not get glued together. 

Now that we have a formal notation and the little observation above, 
we start to play around looking at simple arrows in various categories,
and also at not-so-simple arrows representing standard counterexamples.
 
You notice a few words from your first course on topology:
{\em %injective, surjective, 
$(i)$ connected, $(ii)$ the separation axioms $T_0$ and $T_1$, $(iii)$ dense, $(iv)$ induced (pullback) topology}, and $(v)$ {\em Hausdorff %; surjective and injective on $\pi_0$
} 
are, resp.,
$$\!\!\!\!\!\!\!\!\!\!\!\!\!\!\!\!\!\!\!\!\!\!\!\!\!\!\!\!\!\!\!(i):\ \ \ \ \ \ \ \ \  \ \ \ \ \ \ \ \ \ \ \   \ \ \ \ \ \ \ \  \ \ \   \ \ \ \ \ \ \ \  \ \ \  \ \ X\lra\{\bullet\}\rtt \{\bullet,\bullet\}\lra\{\bullet\}$$
$$(ii):\ \ \ \ \ \ \ \ \ \ \{\bullet\gtrless \star\}\lra\{\bullet\} \rtt X\lra\{\bullet\}\text{\,  and  \,}
\{\bullet< \star\}\lra\{\bullet\} \rtt X\lra\{\bullet\}$$
$$
(iii),(iv):\ \ \ \ \ \ \ \ \
X\lra Y\rtt \{\bullet\}\lra\{\bullet\rightarrow\star\}
\text{\,  and  \,}
X\lra Y\rtt \{\bullet<\star\}\lra\{\bullet\}$$
$$\!\!\!\!\!\!\!\!\!(v):\ \ \ \ \ \ \ \ \  \ \ \ \ \ \ \ \  \ \ \   \ \ \ \ \ \ \ \  \ \   \ \ \ \ \ \ \ \  \ \ \  \  \{\bullet,\bullet'\}\lra X \rtt \{\bullet>\star<\bullet'\}\lra \{\bullet\}$$
See the last two pages for illustrations how to read and draw on the blackboard
these lifting properties  in topology
(here  $\{\bullet<\star\}$, $\{\bullet\gtrless\star\}$, ... denote finite preorders, or, equivalently, 
finite categories with at most one arrow between any two objects, or finite topological spaces
on their elements or objects,
where a subset is closed iff it is downward closed (that is, together with each element, it contains all the smaller elements).
Thus $\{\bullet<\star\}$, $\{\bullet\gtrless\star\}$ and  $\{\bullet>\star<\bullet'\}\lra \{\bullet\}$ 
denote the connected spaces with only one open point $\bullet$,  with no open points, and with two open points
$\bullet,\bullet'$ and a closed point $\star$.
Line $(v)$ is to be interpreted somewhat differently: we consider {\em all} 
the arrows  of form $\{\bullet,\bullet'\}\lra X$).

We mentioned that the lifting property can be seen as a kind of negation. Confusingly, there are {\em two} negations, depending on whether the morphism appears on the left or right
side of the square, that are quite different: for example, both the pullback topology and the separation axiom $T_1$ are 
negations of the same morphism, and the same goes for injectivity and injectivity on $\pi_0$ (see Figure~1(c,d)).

Now consider the standard example of something non-compact: the open covering 
 $\Bbb R= \bigcup\limits_{n\in\Bbb N} \{\,x\,:\, -n<x<n \,\}$ of the real line by infinitely many increasing intervals.
A related arrow in the category of topological spaces is  $\bigsqcup\limits_{n\in\Bbb N} \{\,x\,:\, -n<x<n \}\,\,\lra \Bbb R$. 
Does the lifting property relative to that arrow define compactness? Not quite but almost: 
$$\{\}\lra X \rtt \bigsqcup\limits_{n\in\Bbb N} \{x\,:\, -n<x<n \,\}\,\,\lra \Bbb R$$
reads,  for $X$ connected,  as 
``Every continuous real-valued function on $X$ is bounded", which is an early characterisation of compactness taught in a first course on analysis.
%\vskip 4pt

\,A category theorist would rewrite (**)${}_{\rtt}$ as \bi\item[]
  $(**)_{\text{mono}}\ \ \ \ \ \ \ \ \ \ \ \ \ \ \ \  
  \bullet \!\vee \bullet\lra \bullet \, \rtt \, X\lra Y$
  \ei
denoting by $\vee$ and $\bullet \vee \bullet\lra \bullet$ the coproduct and  the codiagonal morphism, resp., and then
read it as follows: given two morphisms $\bullet\xra{\rm left} X$ and $\bullet\xra{\rm right} X$, 
if the compositions 
$\bullet\xra{\rm left}X\lra Y\ =\ \bullet\xra{\rm right} X\lra Y$ 
are equal (both to $\bullet\xra{\rm down} Y$),
then 
$\bullet\xra{\rm left}X\ =\ \bullet\xra{\rm right} X$ 
are equal (both to $\bullet\xra{\rm down} X$). 
Naturally her first assumption would be that $\bullet$ denotes an {\em arbitrary} object,
as that spares the extra effort needed to invent the  axioms particular to  
the category of sets (or topological spaces) that capture that $\bullet$ denotes
a single element, i.e.~allow to treat $\bullet$ as a single element. 
(A logician understands  ``arbitrary"  as ``we do not know'', ``make
no assumptions'', and that is how formal derivation systems treat ``arbitrary" 
objects.) Thus she would read (**)${}_{\rtt}$ as  the usual category theoretic definition of a monomorphism. 
Note this reading doesn't need that the underlying category has coproducts: a category theorist
would think of working inside a larger category with formally added coproducts $\bullet\vee\bullet$,
and a logician would think of working inside a formal derivation system where ``\,$\bullet$\,'' is 
either a built-in or ``a new variable" symbol,
and ``$\,\bullet \vee \bullet\lra \bullet\,$'' (or ``\,$ \{\bullet,\bullet\}\lra \{\bullet\}$\,'') is (part of) 
a well-formed term or formula.

And of course, nothing prevents a category theorist to make a dual diagram
 \bi\item[]
   $(**)_{\mathrm{epi}}\ \ \ \ \ \ \ \ \ \ \ \ \ \ \ \  
X\lra Y\rtt\, \bullet\lra \bullet \times  \bullet,\ \ \ \ \bullet{\text{ runs through all the objects}} $
\ei
and read it as: 
\bi\item[]$X\lra Y \xra{\rm left}\bullet\ =\ X\lra Y \xra{\rm right}\bullet$ implies $ Y \xra{\rm left}\bullet\ =\ Y \xra{\rm right}\bullet$\ei
which is the definition of an epimorphism.

\section{Speculations.}
Does your brain (or your kitten's) have the lifting property built-in?
Note [G0] suggests a broader and more flexible context making 
contemplating an experiment possible. Namely, some standard 
arguments in point-set topology are computations with 
category-theoretic (not always) commutative diagrams of preorders, in the same way 
that lifting properties define injection and 
surjection. In that approach, the lifting property is viewed as a rule to add a new arrow,
a computational recipe to modify diagrams. 
Moreover, this and some other computations coming from standard theorems in point-set topology
do not involve automorphisms
and thus may perhaps be decidable by an algorithm [GLS] or a modification thereof.
Can one find an experiment  
to check whether humans {\em subconsciously} use diagram chasing to reason about topology? 
Does it appear implicitly in old original papers and books on point-set topology? 
Would teaching diagram chasing hinder or aid development of topological intution
in a first course of topology? Say if one defines connected, dense, Hausdorff et al
via the lifting properties (i-v)?

Is diagram chasing with preorders too complex to have evolved? Perhaps; but note the self-similarity:
preorders are categories as well, with the property that there is at most
one arrow between any two objects; in fact sometimes these categories are 
thought of as 0-categories. So essentially your computations are in 
the category of (finite 0-) categories.

Is it universal enough? Diagram chasing and point-set topology, arguably 
a formalisation of ``nearness'', 
is used as a matter of course in many arguments in mathematics. 
%FIXME: Rather, I want to point out that 
%in many arguments in mathematics are intertwined, mixed in 
%in the (more complicated etc) fabric of the argument. 

Finally, isn't it all a bit too obvious? %Does everybody know these properties or finds them obvious? 
Curiously, in my experience it's a party topic people often get stuck on.
If asked, few if any can define a surjective or an injective map without words,
by a diagram, 
or as a lifting property, even if given the opening sentence of the previous section
as a hint. 
No textbooks seem to bother to mention these reformulations (why?). 
An early version of [GH-I] states (*)${}_{\rtt}$ and (**)${}_{\rtt}$ 
as the simplest examples of lifting properties we were able to think up;
these examples were removed while preparing for publication.

No effort has been made to provide a complete bibliography; 
the author shall happily include any references
suggested by readers in the next version [G].

\subsubsection*{Acknowledgments and historical remarks} It seems embarrassing to thank anyone for ideas so trivial, and 
we do that in the form of historical remarks. 
Ideas here have greatly influenced by extensive discussions with Grigori Mints, Martin Bays, and, later, with 
Alexander Luzgarev and Vladimir Sosnilo. At an early stage Ksenia Kuznetsova
helped to realise an earlier reformulation of compactness was inadequate and that labels on arrows are necessary
to formalise topological arguments. 
``A category theorist [that rewrote] (**) as" the usual category theoretic definition of a  monomorphism,
is Vladimir Sosnilo. %made on 22 Aug 2014 
Exposition has been polished in the numerous conversations with students at St.Peterburg 
and Yaroslavl'2014 summer school.

Reformulations (*)${}_{\rtt}$ and (**)${}_{\rtt}$ of surjectivity and injectivity, as well as connectedness
and (not quite) compactness, 
appeared in early drafts of a paper [GH-I] with Assaf Hasson
 as
trivial and somewhat curious examples of a lifting property 
but were removed during preparation for publication. After (**)${}_{\rtt}$ 
%a similar reformulation of injectivity as a lifting property $\{x,y\}\lra \{x=y\} \rtt X\lra Y$ 
came up in a conversation with Misha Gromov the author decided to try to
think seriously about such lifting properties, and in fact gave talks at logic seminars in 2012 at Lviv
and in 2013 at Munster and Freiburg, and a seminar in 2014 at St.Petersburg.
At a certain point the author realised that possibly a number of simple
arguments in point-set topology may become diagram chasing computations with finite topological spaces, and Grigori Mints
insisted these observations be written. 
 Ideas of [ErgB] influenced this paper (and [GH-I] as well),
 and particularly our computational approach to category theory.
Alexandre Borovik suggested to write a note for [DMG] 
explaining the observation that 
`some
 of human's ``natural proofs'' are expressions
  of lifting properties as applied  to ``simplest counterexample"'.

I thank Yuri Manin for several discussions motivated by [GH-I].

I thank Kurt Goedel Research Center, Vienna, and particulary 
 Sy David Friedman, Jakob Kellner, Lyubomyr Zdomskyy, 
and Chebyshev Laboratory, St.Petersburg, for hospitality.

I thank Martin Bays, Alexandre Luzgarev and Vladimir Sosnilo for proofreading
which have greatly improved the paper. 

 I  wish to express my deep thanks to Grigori Mints, to whose memory this paper is dedicated...

\includepdf[pages=1]{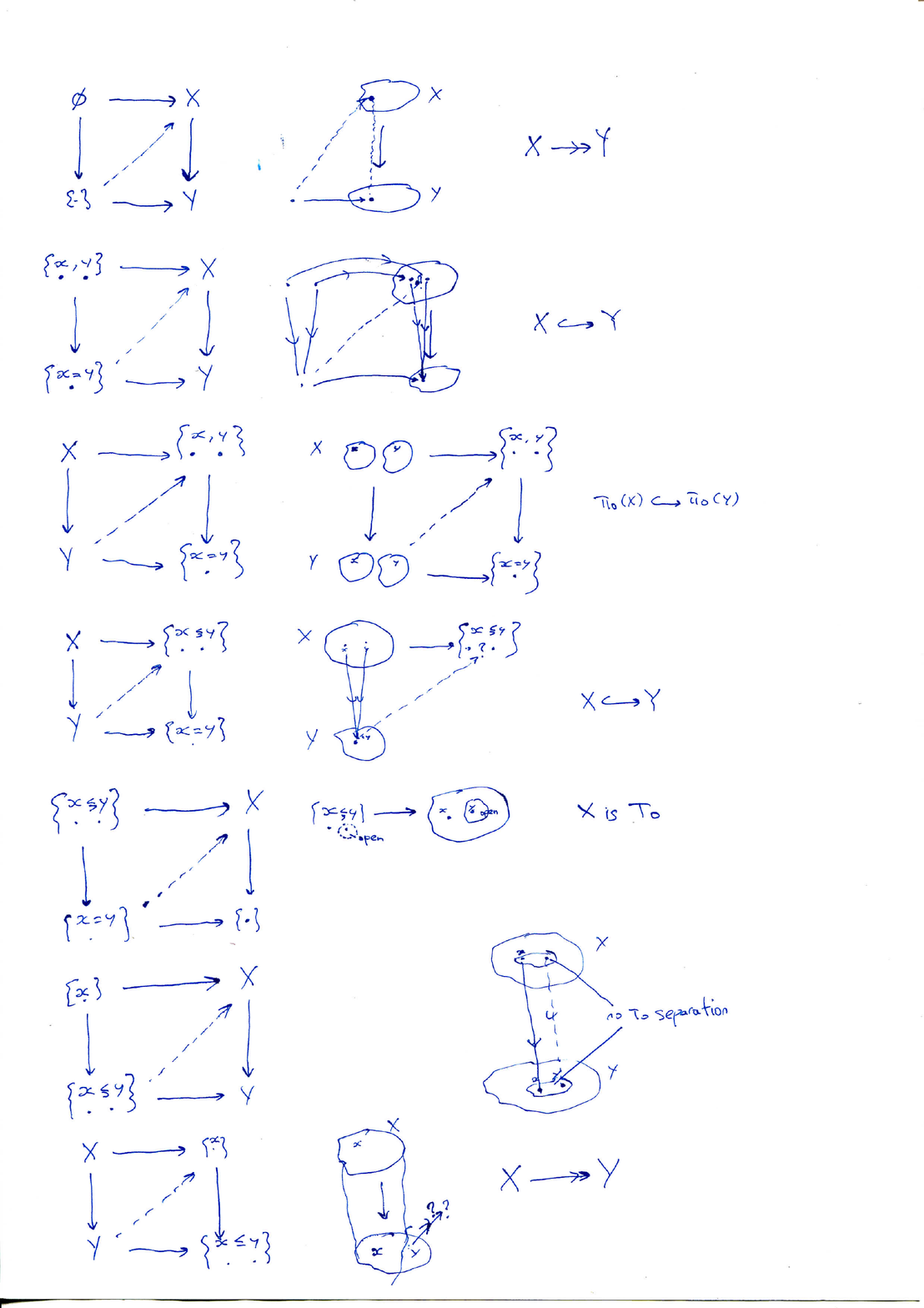}
\includepdf[pages=1]{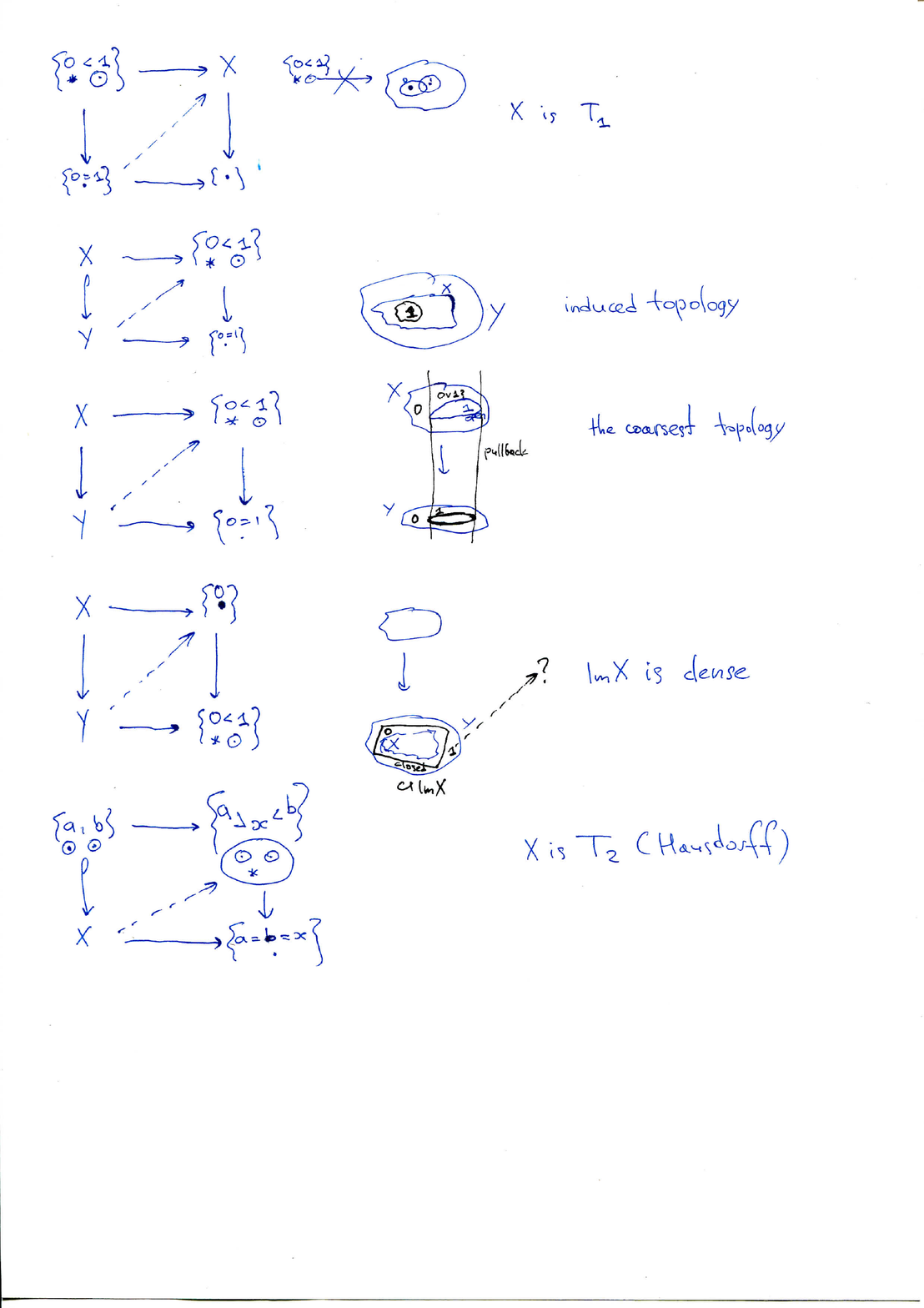}
\end{document}